\documentclass[12pt,aps,prb,preprint,showkeys]{revtex4} 
\usepackage{amsmath} 
\usepackage{amsfonts} 
\usepackage{amssymb}
\usepackage{graphicx}
  
\begin{document}

\title{Interception in differential pursuit/evasion games}

\author{J. A. MORGAN}
\affiliation{The Aerospace Corporation, P. O. Box 92957 \\Los Angeles, CA 90009, 
United States of America \\
\email{john.a.morgan@aero.org}
}

\begin{abstract}

%A recipe for seafood gumbo follows.  

A qualitative criterion for a pursuer to intercept a target in a class of differential games is obtained in terms of 
\emph{future cones}: Topological cones that contain all attainable trajectories of target or interceptor originating from an initial position.  An interception solution exists after some initial time iff  the future cone of the target lies within the future cone of the interceptor.  The 
%condition that gives the existence of the 
solution may be regarded as a kind of Nash equillibrium. This result is applied to two examples:

\noindent 1. The game of Two Cars: The future cone condition is shown to be equivalent to conditions for interception obtained by Cockayne.~\cite{ref2}

\noindent 2. Satellite warfare:  The future cone for a spacecraft or direct-ascent antisatellite weapon (ASAT) maneuvering in a central gravitational field is obtained and is shown to equal that for a spacecraft which maneuvers solely by means of a single velocity change at the cone vertex.

The latter result is illustrated with an analysis of the January 2007 interception of the FengYun-1C spacecraft.

\end{abstract}

\maketitle

\section{Introduction} \label{intro}

A variety of pursuit/evasion problems may  be treated by the methods of differential game
 theory.~\cite{I1965,F1971}  These range from idealized problems of an illustrative nature, such as the Homicidal Chauffeur or Two Cars games,~\cite{I1965,ref2,ref3} to detailed studies of optimal strategies in air-to-air combat.~\cite{ref4,ref5,ref6} In a previous paper~\cite{Morgan2010}, a qualitative criterion for interception or capturability was devised for a class of differential games of kind suggested by Pontryagin.~\cite{P1964,P1971,P1974,P1981} The criterion is a simple one:  
 Let the \emph{future cone} $K^{+}$ be the set of all attainable trajectories available to a player subsequent to some initial time.  Interception of the target by the interceptor is \emph{guaranteed} to be possible if 
\begin{equation}
K^{+}_{target} \subset K^{+}_{interceptor}
\end{equation}
subject to additional assumptions given below.  The condition takes the form of a dominant strategy Nash equillibrium, in that, no matter how the target maneuvers, the interceptor can always maneuver to intercept it at some point to the future of both.

Section \ref{newtheorem} provides a new proof of the main result of Ref.  \onlinecite{Morgan2010}, based on the Lefschetz fixed-point theorem, which relaxes assumptions regarding the convexity of subsets of the future cones made in  Ref.  \onlinecite{Morgan2010}.  In Section \ref{futurecones}, the criterion is applied to the classic problem of the game of Two Cars, and to the problem of two spacecraft maneuvering in a central gravitational field.  

This latter problem is applicable to the study of a direct-ascent antisatellite weapon (ASAT) engaging a spacecraft in low Earth orbit.  An account of the January 2007 interception of the FengYun-1C meteorological satellite by a Chinese ASAT appears in Section \ref{FY1C} as a worked example.

Following a discussion of the results in Section \ref{wottdoesitmean}, a final Section \ref{byby} presents conclusions.

\section{Conditions for guaranteed interception}  \label{newtheorem}

This paper is concerned with differential games that describe the pursuit of a target with position $y(t)$ at time $t$, by an interceptor whose position at time $t$ is $x(t)$.  Both target and interceptor can maneuver freely and autonomously, subject to constraints.  The evolution of $x$ and $y$ is of the form
 \begin{equation} 
 \frac{dx}{dt}=F(x,u) \label{eq:xeqn}
 \end{equation} 
 and
  \begin{equation} 
  \frac{dy}{dt}=G(y,w). \label{eq:yeqn} 
  \end{equation} 
 where $F$ and $G$ are assumed to be bounded analytic functions, and $u=u(t)$ and $w=w(t)$ are piecewise analytic controls. 

\subsection{The Future Cone} \label{futurecone}

We recall the definition of the future cone
of a maneuvering player as given in Ref. \onlinecite{Morgan2010}, and review its properties. The net effect of propulsive forces and of external forces such as gravity or aerodynamic drag cause the position $x(t)$ of the interceptor (respectively, target $y(t)$) to evolve as it maneuvers.  The evolution will be continuous, but not necessarily differentiable.  
In addition, target and interceptor may be subject to physical limitations on their peak acceleration or total velocity change $\Delta v$. Finally, if target and interceptor maneuver freely in $\mathbf{R}^{3}$ during a finite interval of time, they may be expected to interact within a compact subset of $\mathbf{R}^{4}$.  

At an initial time $t_{0}$, call the position of the interceptor  $x(t_{0})$.  
The set of all possible histories for $x(t)$ originating at $x(t_{0})$ comprises a 
topological cone in $\mathbf{R}^{4}$ with vertex $x(t_{0})$. We call the set of points subsequently accessible to the interceptor in the time interval $(t_{1},t_{2})$, with $t_{0} \le t_{1} < t_{2}$ the \emph{future cone} of $x(t_{0})$, written 
$K_{x}^{+}(t_{1}:t_{2};x(t_{0}))$.\footnote{If the future cone is a closed set, it is identical to the \emph{attainable set} introduced by Roxin in Ref. \onlinecite{Roxin1969} for a different purpose.} 

Following Ref. \onlinecite{Morgan2010}, we assume that the subset of either cone lying to the future of its vertex is a manifold with compact closure possessing a timelike foliation. We shall
denote by  $K^{+}_{x}(t;x(t_{0}))$ a leaf of the foliation corresponding to a time $t$.
In particular, there is no assumption that future cones of either target or interceptor, or their leaves, are necessarily convex.
 
\subsection{Guaranteed Interception}

The properties of optimal solutions to (\ref{eq:xeqn}) and (\ref{eq:yeqn}) leading to 
interception has been well-studied since Ref.  \onlinecite{P1964}; \emph{vide.} Refs. \onlinecite{P1971,P1974,P1981}.
A sizeable body of  literature is devoted to the closely related Homicidal 
Chauffeur and Two Cars games introduced by Isaacs  in Ref \onlinecite{I1965}; \emph{e. g.} Refs. \onlinecite{ref1,ref5,ref4,ref3,ref6,ref2}. These studies variously treat games of degree or of kind, but all restrict the form of the game to facilitate its  analysis, by means such as limitation to a
linear game, or to fixed velocity ratios for the players, or to piecewise constant radii of 
curvature. 

We work here with (\ref{eq:xeqn}) and (\ref{eq:yeqn}) in fairly general form, but restrict attention
to the more limited goal of finding qualitative  conditions 
under which we may be confident that the interceptor can force an interception, optimally or no. The winning strategy for the interceptor is to 
choose a trajectory that intercepts the target at some time to the future of  $t_{0}$.
This choice amounts to a mapping
$\mathbf{\Phi}:K^{+}_{x} \rightarrow K^{+}_{x}$ from the set of all trajectories available to the interceptor 
to its desired \emph{actual} trajectory.
A \emph{guaranteed intercept} will
be said to exist when the pursuer always has available to it a strategy leading to interception, no matter how the target maneuvers 
within its future cone. 
%We shall suppose that interception occurs with high probability if the interceptor position $x(t_{i})$ equals the target position $y(t_{i})$ at the time of interception $t_{i}$. This criterion serves to allow a simple formulation of the game while in some wise accommodating the circumstance that in actuality an interceptor can only navigate to a stipulated position at any given time within some nonzero error.

\subsection{New proof based on fixed-point theorem for correspondences on nonconvex domains}

We use the Eilenberg-Montgomery theorem~\cite{EM1946,Tefatsion} to prove a version of Theorem 3.3 from 
Ref.  \onlinecite{Morgan2010} giving the conditions for the existence of a guaranteed intercept:
\newline
\,

\noindent Theorem 1:  \emph{
Let 
$t_{\alpha} < t_{0} < t_{1} < t_{\omega}$. Then, 
a necessary and sufficient condition for the existence of a guaranteed intercept at time of interception 
$t_{i} \in [t_{0},t_{1}]$ is}
\begin{equation}
K^{+}_{y}(t_{0}:t_{1};y(t_{0})) 
\subset K^{+}_{x}(t_{\alpha}:t_{\omega};x(t_{\alpha})) \label{eq:conecond}
\end{equation}

\emph{Proof}: Sufficiency: Suppose that, of all the possible trajectories 
$\subset K^{+}_{y}(t_{i};y(t_{0}))$, the actual trajectory of the target is $y^{*}(t)$. The mapping 
from  $K^{+}_{x}(t_{\alpha}:t_{\omega};x(t_{\alpha}))$  into
$K^{+}_{x}(t_{\alpha}:t_{\omega};x(t_{\alpha}))  \cap K^{+}_{y}(t_{0}:t_{1};y(t_{0}))=
K^{+}_{y}(t_{0}:t_{1};y(t_{0}))$ is given by the correspondence
$\mathbf{\Phi}(x) =  \{y^{*}(t), t \in [t_{0},t_{1}] \}$, which is continuous.  Its value is homeomorphic to a one-simplex in 
%$\mathbf{\Phi}(x):  K^{+}_{x}(t_{\alpha}:t_{\omega};x(t_{\alpha})) \rightarrow  \{y^{*}(t), t \in [t_{0},t_{1}] \}$, which is continuous.  Its value is homeomorphic to a one-simplex in 
$\mathbf{R}^{4}$, and is thus acyclic. The cones  $K^{+}_{x}(t_{\alpha}:t_{\omega};x(t_{\alpha}))$ and $K^{+}_{y}(t_{0}:t_{1};y(t_{0}))$ are differentiable manifolds by construction, and thus
triangulable.~\cite{Cairns46}  Each may therefore be regarded as the polyhedron of a simply connected simplicial complex and is, in consequence, both an acyclic set~\cite{Spanier66} and an absolute neighborhood retract.~\cite{GranasDj03}
The conditions for the Eilenberg-Montgomery theorem are satisfied:  A fixed point of the mapping $\mathbf{\Phi}(x)$
\begin{equation}
x^{*}(t_{i}) \in \{y^{*}(t), \forall t \in [t_{0},t_{1}] \}
\end{equation}
exists for 
$t_{0} < t_{i} < t_{1}$.
Combining this result with the tautological fixed point 
$y^{*}(t_{i}) \in K^{+}_{y}(t_{0}:t_{1};y(t_{0}))$ resulting from the target's ability to maneuver
freely, we find
\begin{equation}  
\left( \begin{array}{c} x^{*}(t_{i})\\ y^{*}(t_{i}) \end{array} \right) \in 
\left( \begin{array}{c} \{y^{*}(t), \forall t \in [t_{0},t_{1}] \} \\  
K^{+}_{y}(t_{0}:t_{1};y(t_{0})) \end{array} \right). 
\label{eq:Nash2}
\end{equation}

The proof of the necessary condition is identical to that in Ref. \onlinecite{Morgan2010} but is included here to make the treatment self-contained.
We begin with a condition for the relation between leaves of 
the two cones at a single time. A necessary condition for the existence of a guaranteed intercept is that, at the time of intercept $t_{i}$,
\begin{equation}
K^{+}_{y}(t_{i};y(t_{1})) \subseteq K^{+}_{x}(t_{i};x(t_{0})). \label{eq:leafcond}
\end{equation}
for $t_{0}, t_{1} < t_{i}$:
Suppose a guaranteed intercept exists at 
time $t_{i}$.  Then, every point $y$ in $K^{+}_{y}(t_{i};y(t_{1}))$ must 
coincide with some point $x$ in 
$K^{+}_{x}(t_{i};x(t_{0}))$ in order that $\parallel x-y \parallel = 0$ for at least one pair of 
values of $x$ and $y$. Were (\ref{eq:leafcond}) false, there
would be some portion of $K^{+}_{y}(t_{i})$ that lay outside the attainable set of interceptor
positions at that time.  Thus there would be a subset of $K^{+}_{y}(t_{i};y(t_{1}))$ for which
$\parallel x-y \parallel > 0, \forall x \in K^{+}_{x}(t_{i};x(t_{0}))$.

The necessary condition for the entire cone $K^{+}_{y}(t_{0}:t_{1};y(t_{1}))$ is obtained by transfinite induction \citep{K1955}.  Let 
$t_{\alpha} < t_{0}  < t_{\beta} < t_{\gamma} < t_{i} < t_{1} < t_{\omega}$ 
and take $t_{i}-t$ as an ordinal.  We prove the result for 
$K^{+}_{y}(t_{\beta}:t_{1};y(t_{0})) \subset K^{+}_{y}(t_{0}:t_{1};y(t_{0}))$ and extend to
the full set $K^{+}_{y}(t_{0}:t_{1};y(t_{0}))$ at the end.

We begin by showing the necessary condition holds at late times.  
Suppose that a guaranteed intercept exists at time $t_{i}$ for $t_{\gamma}, t_{i}$ 
within any neighborhood of $t_{1}$.
As $t_{\gamma} \rightarrow t_{1}$,
\begin{equation} 
K^{+}_{y}(t_{\gamma}:t_{1};y(t_{0})) \rightarrow K^{+}_{y}(t_{1};y(t_{0})).
\end{equation}  
By the condition for a single leaf, it follows that
\begin{equation} 
K^{+}_{y}(t_{1};y(t_{0})) \subseteq K^{+}_{x}(t_{1};x(t_{\gamma}))
\subset K^{+}_{x}(t_{\alpha}:t_{\omega};x(t_{\alpha})) 
\end{equation}
  
Next, suppose that at least one guaranteed intercept opportunity
exists for time $t_{i}$ between $t_{\gamma}$ and $t_{1}$.  By the inductive hypothesis,
\begin{equation}
K^{+}_{y}(t_{\gamma}:t_{1};y(t_{0})) 
\subset K^{+}_{x}(t_{\alpha}:t_{\omega};x(t_{\alpha})) 
\end{equation}
We wish to examine the prospects at 
an earlier time $t_{\beta}$.  Consider the sets 
$K^{+}_{y}(t_{\beta}:t_{\gamma};y(t_{0}))$ and $K^{+}_{x}(t_{\beta}:t_{\gamma};x(t_{\alpha}))$:  
\begin{equation}
K^{+}_{y}(t_{\beta}:t_{1};y(t_{0}))=K^{+}_{y}(t_{\beta}:t_{\gamma};y(t_{0})) \cup
K^{+}_{y}(t_{\gamma}:t_{1};y(t_{0})) \label{eq:eqn11}
\end{equation}
and similarly for $K^{+}_{x}$.  But if a guaranteed intercept is to be possible  
$\forall t \in (t_{\beta},t_{\gamma})$, at no
time $t$ in $(t_{\beta},t_{\gamma})$ can it be that
\begin{equation} 
K^{+}_{y}(t;y(t_{0})) \not \subseteq K^{+}_{x}(t;x(t_{\alpha})),
\end{equation}
by the single-leaf condition.  Letting $t_{\beta} \rightarrow t_{0}$ in (\ref{eq:eqn11}) and recalling $t_{\alpha} < t_{0}$, we have (\ref{eq:conecond}).
$\,\square$

One may regard the union of leaves $K^{+}_{y}(t;y(t_{0}))$ or  $K^{+}_{x}(t;x(t_{\alpha}))$ comprising the future cones in (\ref{eq:conecond}) as subsets of either 
$\mathbf{R}^{4}$ or $\mathbf{R}^{3}$; the latter case amounts to a projection $\mathbf{R}^{4} \rightarrow \mathbf{R}^{3}$. The projection into $\mathbf{R}^{3}$ lends itself
to simple graphical presentation.

As noted in Ref. \onlinecite{Morgan2010}, the result (\ref{eq:Nash2}) may be interpreted as a Nash equillibrium.~\cite{N1950} The interceptor strategy in (\ref{eq:Nash2}) is a strictly dominant one.~\cite{N1950a} That is, while the strategy of the target allows it to move anywhere within its future cone, the interceptor always has available to it the strategy which places it at some future position of the target. 
%If we stipulate $t_{0}$ in $K^{+}_{y}(t;y(t_{0}))$, 
On the other hand, no alternative to the target's strategy given in (\ref{eq:Nash2}) will increase its chances of survival. We note that, having stipulated $t_{0}$ in $K^{+}_{y}(t;y(t_{0}))$, any such alternative will necessarily constrict the volume available for maneuvering by the target, with (we may suppose) the effect of worsening its prospects. 
%We argue that the target strategy may likewise be considered strictly dominant. As mentioned earlier, we assume that the existence of a guaranteed interception guarantees only an opportunity to intercept; there is a small probability that the interceptor will fail to strike the target at $y^{*}(t_{i})$. If we stipulate $t_{0}$ in $K^{+}_{y}(t;y(t_{0}))$, any alternative to the target strategy given in (\ref{eq:Nash2}) will reduce the integrated probability of a failed interception. 
%If we accept this consideration,
If we argue in this way,
the strategy in (\ref{eq:Nash2}) is preferable to any other available to the target, and thus strictly dominant. Both target and interceptor then have 
available to them 
a strictly dominant, hence optimal, strategy. We may conclude (\ref{eq:Nash2}) is the unique Nash equillibrium for this problem.

\section{Calculation of Future cones for Selected Problems} \label{futurecones}

\subsection{Two Cars} \label{2cars}

In this game, introduced by Isaacs~\cite{I1965}, two cars maneuver in a plane. We have a target Car 2 with position $r_{2} \equiv (x_{2}(t),y_{2}(t))$ at time $t$, pursued by an interceptor Car 1 whose position at time $t$ is $r_{1} \equiv  (x_{1}(t),y_{1}(t))$.  
Players move with constant velocities $v_{1}$ and 
$v_{2}$, respectively, and maneuver exclusively by steering. The motion of Car 1 in the plane is given by 
\begin{equation}
\frac{dx_{1}}{dt}=v_{1} \, sin(\theta_{1}(t)) \label{eq:xeqn2}
\end{equation}
and
\begin{equation}
\frac{dy_{1}}{dt}=v_{1} \, cos(\theta_{1} (t)). \label{eq:yeqn2}
\end{equation}
The control law for $\theta_{1} (t)$ is subject to the constraint
\begin{equation}
| \dot{\theta_{1}} | < \frac{v_{1}}{R_{1}}. \label{eq:thetadotcond}
\end{equation}
Identical laws govern the motion of Car 2.
  
Cockayne~\cite{ref2} showed that  sufficient and necessary conditions for Car $1$ to intercept Car $2$ are
\begin{equation}
 v_{1} > v_{2}. \label{eq:vcond}
\end{equation}
and 
\begin{equation}
 \frac{v_{1}^{2}}{R_{1}} \ge  \frac{v_{2}^{2}}{R_{2}}. \label{eq:acond}
\end{equation} 
We show that the condition for guaranteed interception given by applying Theorem 1 
\begin{equation}
K^{+}_{2} \subset K^{+}_{1}
\end{equation}
to the game of two cars is equivalent to the paired conditions (\ref{eq:vcond})-(\ref{eq:acond}) found by Cockayne.  It is convenient to assume that the initial velocity vector of either car may point in any direction.  This assumption ensures that the future cones are simply connected.
 
To see equivalence of the sufficient conditions, we adapt Isaacs'~\cite{I1965} construction of trajectories for Car 1 that result in interception of Car 2. 
Amongst the admissible trajectories $\in K^{+}_{1}$ must be those produced by the 
Method of the Explicit Policy.  The account of the method as applied to this problem given by Isaacs~\cite{I1965} can hardly be bettered (P is the pursuer, Car 1 and E, the evader, is Car 2):

 \emph{If P has the higher speed and at least as favorable a curvature restriction as E, capture can be attained.  For P can first go to E's starting point and then follow his track.}

It is clear that, if $K^{+}_{2}
\subset K^{+}_{1}$, Car 1 can accomplish this for every trajectory of Car 2 $\in  K^{+}_{2}$.
Let $K^{+}_{XP}$ be the cone comprised of all trajectories constructed by Explicit Policy:
\begin{equation}
 K^{+}_{1} \supset K^{+}_{XP}=K^{+}_{2}.
\end{equation}
The resulting sufficient conditions required for the motion of Car 1 are two~\cite{I1965}: 

\noindent 1.  $v_{1} > v_{2}$:
Assume
\begin{equation}
K^{+}_{2}(t_{2}:t_{f};r_{2}(t_{1}))
\subset K^{+}_{1}(t_{1}:t_{f};r_{1}(t_{0})). \label{eq:hypothesis}
\end{equation}
with 
\begin{equation}
t_{0} < t_{1} < t_{2} \label{eq:headstart}
\end{equation}
and $t_{f} > t_{2}$ large, but otherwise arbitrary. In (\ref{eq:headstart}), one may suppose that the difference between $ t_{2}$ and $ t_{1}$ includes an overestimate of the time for Car 1 to come about, should its initital velocity not be parallel to that of the target. 

Within $K^{+}_{1}(t_{1}:t_{f};\mathbf{r}_{1}(t_{0}))$ consider the future cone for Car 1 originating at 
$\mathbf{r}_{2}(t_{1})$.  
We wish to force an interception for $t > t_{2}$. By hypothesis,
\begin{equation}
K^{+}_{2}(t_{2}:t;r_{2}(t_{1}))
\subset K^{+}_{1}(t_{2}:t;r_{1}(t_{2})). \label{eq:hypothesis2}
\end{equation}
Car 1's position at time $t$ may be written formally as
\begin{equation}
x_{1}(t)=v_{1} \, \int_{t_{0}}^{t}sin(\theta_{1}(s))  \,  ds   \label{eq:xmotioneqn}
%+ x_{1}(t_{0})
\end{equation}
and
\begin{equation}
y_{1}(t)=v_{1} \, \int_{t_{0}}^{t}cos(\theta_{1}(s))  \,  ds . \label{eq:ymotioneqn}
%+ y_{1}(t_{0})
\end{equation}
The distance traveled from the origin at time t is
\begin{eqnarray}
\, & r_{1}(t)= v_{1} \,\sqrt{ x_{1}^2+y_{1}^{2}} & \,  \label{eq:rmotioneqn} \\  
 \, & =v_{1} \, \sqrt{\left ( \int_{t_{0}}^{t}sin(\theta_{1}(s)) \,  ds  \right )^2+
 \left (\, \int_{t_{0}}^{t}cos(\theta_{1}(s)) \,  ds  \right )^2}.  & \,
\end{eqnarray}
The Schwartz inequality gives
\begin{equation}
\left (  \int_{t_{0}}^{t}sin(\theta_{1}(s)) \,  ds   \right )^{2} \le  
\int_{t_{0}}^{t}sin^{2}(\theta_{1}(s)) \,  ds \,  \int_{t_{0}}^{t} \,  ds 
\end{equation}
and similarly for the cosine. We thus find
\begin{equation}
r_{1}(t) t_{1} \le   v_{1}  \,  \sqrt{ \int_{t_{0}}^{t} \,  ds} \, 
\sqrt{ \int_{t_{0}}^{t} (sin^{2}(\theta_{1}(s)) +
cos^{2}(\theta_{1}(s))) \,  ds}
\end{equation}
or
\begin{equation}
r_{1}(t_{1}) \le   v_{1}  \, (t_{1}-t_{0}). \label{eq:rcondition}
\end{equation}

By (\ref{eq:rcondition}), the leaf $K^{+}_{2}(t_{2})$ of the target at a time $t > t_{2}$ is contained within a circle of radius $r_{2}$, and the leaf $K^{+}_{1}(t_{2})$ within a circle of radius $r_{1}$.  From (\ref{eq:hypothesis2}),
\begin{equation}
r_{2}=v_{2}(t-t_{2}) < r_{1}= v_{1}(t-t_{2})
\end{equation}
or
\begin{equation}
 v_{1} > v_{2}. \label{eq:vcond2}
\end{equation}

\noindent 2. $ \frac{v_{1}^{2}}{R_{1}} \ge  \frac{v_{2}^{2}}{R_{2}}$:
The Method of the Explicit Policy requires
\begin{equation}
| \dot{\theta_{2}} | <  | \dot{\theta_{1}} | \label{eq:thetadotcond2}
\end{equation}
From (\ref{eq:xeqn2}) and (\ref{eq:yeqn2}),
\begin{equation}
\frac{dx^{2}_{1}}{dt^{2}}=v_{1} \, cos(\theta_{1}(t)) \, \dot{\theta_{1}},  \label{eq:x2eqn2}
\end{equation}
\begin{equation}
\frac{dy^{2}_{1}}{dt^{2}}=-v_{1} \, sin(\theta_{1} (t))\dot{\theta_{1}}, \label{eq:y2eqn2}
\end{equation}
and for Car 2, likewise, whence  the acceleration
\begin{eqnarray}
\, & a= \sqrt{ (v_{1} \, cos(\theta_{1}(t)) \, \dot{\theta_{1}})^2
+(v_{1} \, sin(\theta_{1} (t))\dot{\theta_{1}})^2} & \,  \\  
 \, & = v_{1} | \dot{\theta_{1}} | & \, \\
 \, & \le \frac{v_{1}^{2}}{R_{1}},  & \,  \label{eq:aeqn}
\end{eqnarray}
by (\ref{eq:thetadotcond}). Note that $\mathbf{a} \cdot \mathbf{v}=0$ at all times.
Equations (\ref{eq:thetadotcond2}) and (\ref{eq:aeqn})  immediately give
\begin{equation}
 \frac{v_{1}^{2}}{R_{1}} \ge  \frac{v_{2}^{2}}{R_{2}}. \label{eq:acond2}
\end{equation} 
The sufficient conditions (\ref{eq:vcond2}) and (\ref{eq:acond2}) are  those found by Cockayne.

We may see the equivalence of the necessary condition for guaranteed interception
\begin{equation}
K^{+}_{2} \subset K^{+}_{1} \label{eq:conecond2}
\end{equation}
and (\ref{eq:vcond})-(\ref{eq:acond}) by considering the contrapositive:
 If the conditions be violated, then by Cockayne's analysis, there are points the Car 2 can reach that Car 1 cannot (at certain constructable times).  If that be true, then 
$K^{+}_{2} \not \subset K^{+}_{1}$. 
 
\subsection{Spacecraft maneuvering in a gravitational field} \label{Kepler}

In this section, we calculate the future cone of a spacecraft maneuvering in a gravitational field by  impulsive velocity changes.  The result yields an account of a direct-ascent ASAT engagement with a target satellite.  It is also applicable to the endgame of a co-orbital ASAT 
attack.~\cite{carter1985}

The motion of the spacecraft is given by
\begin{eqnarray}
 \frac{d}{dt} \left ( \begin{array}{c}
 \mathbf{r} \\ \mathbf{v} 
  \end{array} \right ) =
 \left ( \begin{array}{c}
 \mathbf{v} \\ -\mu \frac{\mathbf{r}}{r^{3}} +\mathbf{a}(t)
  \end{array} \right ) \label{eq:eom0}
\end{eqnarray}
where  $\mathbf{r}$ is the spacecraft radius vector, $\mathbf{v}$ is its velocity,  $\mathbf{a}(t)$ its acceleration, and $\mu=G_{Newt.} M_{\oplus}$.  (Because the spacecraft state descriptor includes the velocity, the condition for interception is modified in an obvious way.)
We make use of an exact solution for the system (\ref{eq:eom0}) in the case of ballistic motion in the Earth's gravitational field. The motion of a spacecraft in a bound Keplerian orbit is obtained using Lagrange coefficients $ \mathbf{\Phi} $~\cite{Battin}
\begin{equation}
\left ( \begin{array}{c}
\mathbf{r}(t) \\  \mathbf{v}(t)
\end{array} \right ) =
 \mathbf{\Phi} 
\left ( \begin{array}{c}
\mathbf{r}_{0} \\  \mathbf{v}_{0}
\end{array} \right )
\end{equation} 
with
\begin{equation}
\mathbf{\Phi}=
\left ( \begin{array}{cc}
F(t) & G(t) \\  
F_{t}(t) & G_{t}(t)
\end{array} \right ) \label{eq:Lagrange}
\end{equation} 

It is not possible to express the coefficients in  (\ref{eq:Lagrange}) in closed form as functions of time.  However, it is possible to do so in terms of the corresponding true anomaly $f$.
 If $f_{0}$ is the true anomaly for which
\begin{equation}
\left ( \begin{array}{c}
\mathbf{r}(t(f_{0})) \\  \mathbf{v}(t(f_{0}))
\end{array} \right ) =
\left ( \begin{array}{c}
\mathbf{r}_{0} \\  \mathbf{v}_{0}
\end{array} \right )
\end{equation} 
and
\begin{equation}
\theta \equiv f-f_{0},
\end{equation} 
 we have:
\begin{equation}
F=1-\frac{r}{p} \left ( 1-cos(\theta) \right )
\end{equation}
\begin{equation}
G=\frac{r \, r_{0}}{\sqrt{\mu \, p}} \, sin(\theta)
\end{equation}
\begin{equation}
F_{t}=\frac{\sqrt{\mu}}{r_{0} \, p} \left [ \sigma_{0} (1-cos(\theta) -\sqrt{p} \, sin(\theta) \right ]
\end{equation}
\begin{equation}
G_{t}=1-\frac{r_{0}}{p} \left ( 1-cos(\theta) \right )
\end{equation}
In the foregoing, $a$ is the semimajor axis, $e$ the eccentricity, $r_{0}$ is the initial orbital radius at true anomaly 
$f=f_{0}$, the parameter $p \equiv  \frac{(\mathbf{r} \times \mathbf{v})^2}{\mu}=a(1-e^2)$, and
\begin{equation}
\sigma_{0} \equiv \frac{\mathbf{r_{0}} \cdot \mathbf{v_{0}}}{\sqrt{\mu}}.
\end{equation}

A complete description of the spacecraft motion also requires the eccentric and mean anomalies.  These quantities are related by
\begin{equation}
tan \left (\frac{f}{2} \right )=\sqrt{\frac{(1+e)}{(1-e)}}tan \left (\frac{E}{2} \right ). \label{eq:f2t1}
\end{equation}
(Note that $\frac{f}{2}$ and $\frac{E}{2}$ always lie in the same quadrant.)  
The relation between $E$ and the mean anomaly $M$ given by Kepler's equation:
\begin{equation}
M=E-e \, sin(E). \label{eq:f2t2}
\end{equation}
By use of the mean motion 
\begin{equation}
n=\sqrt{\frac{\mu}{a^{3}}} \label{eq:f2t3}
\end{equation}
(obtained from Kepler's third law) the time $t$ corresponding to true anomaly $f$ is obtained from $M$ and the time of pericenter passage $\tau$ by
\begin{equation}
M=n(t-\tau). \label{eq:f2t4}
\end{equation}
The time interval between two points on a ballistic arc is thus
\begin{equation}
t_{2}-t_{1}=\sqrt{\frac{a^{3}}{\mu}} \, (E_{2}-E_{1}-e \, (sin(E_{2})-sin(E_{1})). \label{eq:deltat}
\end{equation}

We wish to calculate the future cone $K^{+}(t_{1}:t_{2};\mathbf{r}(t_{0}))$ of a spacecraft maneuvering
in the Earth's gravitational field. To that end, we model the motion of an idealized target or interceptor  with a trajectory consisting of piecewise ballistic orbital segments punctuated by shocks in which
\begin{equation}
\mathbf{v}(t^{+})=\mathbf{v}(t^{-})+\Delta \mathbf{v}. \label{eq:shock}
\end{equation}
The acceleration during impulsive maneuvers of actual spacecraft will be bounded and continuous, even when occurring in times much shorter than any other timescale characteristic of the motion. The description of impulsive velocity changes as instantaneous shocks is an idealization adopted for computational convenience. However, in the proof of Theorem 3 below, the acceleration history of the spacecraft is assumed continuous, however violent, in order to satisfy requirements of the Gronwall inequality.  In order to include the limiting case of instantaneous shocks, we may treat the acceleration history as a distribution. The transition in (\ref{eq:shock}) may then be regarded as the limit of a suitable sequence of good functions, and that limit may be taken at the end of other calculations.  

We assume the spacecraft maneuvers by a number (possibly large) of small-impulse shocks, and that the total impulse available to it is limited by
\begin{equation}
\sum_{i=0}^{n-1} \| \Delta \mathbf{v}_{i} \| \le \Delta v_{tot}, \label{eq:many_shox}
\end{equation}
where $ \|  \mathbf{x} \|^{2}=\mathbf{x} \cdot \mathbf{x}$.
The spacecraft trajectory is thus a chain comprised of piecewise elliptical ballistic orbits that is continuous, but not necessarily differentiable at any point corresponding to a shock.  For all values of the true anomaly between shocks $t_{n} < t < t_{n+1}$ the motion is given by the free-fall Lagrange coefficients, in expressions to be given, as a function of
$\theta=f-f^{n}_{0}$.  
In any ballistic interval of the spacecraft motion, $f$ is the true anomaly of the osculating elliptical orbit that coincides with the spacecraft trajectory (in this connection, by "osculating" one should perhaps understand  "liplocked" ).

We require a relation between the true anomaly and time at any point in the course of an engagement.  
A closed-from solution for the time (found by Kepler) may be expressed in terms of the eccentric and mean anomalies.  
Begin with the relations between true and eccentric anomalies within the $n^{th}$ ballistic arc for three points
on the orbit: 
\begin{equation}
tan \left (\frac{f-f^{n}_{0}}{2} \right )=\sqrt{\frac{(1+e)}{(1-e)}}tan \left (\frac{E-E^{n}_{0}}{2} \right ). \label{eq:f2t1a}
\end{equation}
and
\begin{equation}
tan \left (\frac{f^{n}-f^{n}_{0}}{2} \right )=\sqrt{\frac{(1+e)}{(1-e)}}tan \left (\frac{E^{n}-E^{n}_{0}}{2} \right ). \label{eq:f2t1b}
\end{equation}
In the foregoing, $\theta \equiv f-f^{n}_{0}$ (respectively, $E-E_{0}^{0}$) is the difference in true (respectively, eccentric) anomalies between the position $\mathbf{r}$ and the semimajor axis of the $n^{th}$ ballistic arc, while $\theta^{n} \equiv f^{n}-f^{n}_{0}$ 
 ($E_{n}-E^{n}_{0}$) is the difference in true (eccentric) anomalies between the position $\mathbf{r}_{n}$  and 
 semimajor axis.  The time interval between $\mathbf{r}$  and $\mathbf{r_{n}}$ is given by (\ref{eq:deltat}) as
\begin{equation}
t-t_{n}=\sqrt{\frac{a^{3}}{\mu}} \,\left  ((E-E^{n}_{0})-(E^{n}-E_{0}^{0})
-e \, (sin(E-E_{n}^{0})-sin(E_{n}-E_{0}^{0})) \right ). \label{eq:deltat2}
\end{equation}

At a time $t_{n+1}$, let a velocity change $\Delta \mathbf{v}_{n+1}$ be imparted when the spacecraft is at position
 $\mathbf{r}_{n}$.  The initial position and velocity are thus $\mathbf{r}_{n}$ and $\mathbf{v}+\Delta \mathbf{v}_{n+1}$. The semimajor axis of the new ballistic arc is obtained from the  \emph{vis-viva} equation
\begin{equation}
\| \mathbf{v}+\Delta \mathbf{v}_{n+1} \|^2=\mu \left ( \frac{2}{\| \mathbf{r}_{n} \|}-\frac{1}{a_{n+1}} \right).
\end{equation}
The parameter (and hence the eccentricity) is given by 
\begin{equation} 
p=a(1-e^{2})=\frac{\|  \mathbf{r}_{n} \times (\mathbf{v}+\Delta \mathbf{v}_{n+1}) \|^{2}}{\mu},
\end{equation}
and 
\begin{equation}
\sigma_{0} = \frac{\mathbf{r_{n}} \cdot  (\mathbf{v}+\Delta \mathbf{v}_{n+1})}{\sqrt{\mu}}.
\end{equation}
These quantities determine the Lagrange coefficient matrix $\Phi_{n+1}$ for the $n+1^{st}$ ballistic arc.

Call
$K_{n}^{+}(t_{1}:t_{2};\mathbf{r}(t_{0}))$ the set of all trajectories originating from  
$\mathbf{r}(t_{0})$ and subject to $n$ shocks $\{ \Delta \mathbf{v}_{i}, i=1 \ldots n, \}$ respecting (\ref{eq:many_shox}), with the understanding that for $n=1$ the shock occurs at time $t_{0}$.  The
future cone of a spacecraft maneuvering by ballistic motion between shocks is then 
\begin{equation}
K^{+}_{\infty}(t_{1}:t_{2};\mathbf{r}(t_{0})) = \bigcup_{n=1}^{\infty} K_{n}^{+}(t_{1}:t_{2};\mathbf{r}(t_{0})).
\end{equation}

We now prove a result relating the future cone $K^{+}_{\infty}(t_{1}:t_{2};\mathbf{r}(t_{0}))$ of a spacecraft that maneuvers in a gravitational field by free fall punctuated with a countable sequence of shocks
$\{ \Delta \mathbf{v}_{i} \}$, subject to the overall upper limit on their sum (\ref{eq:many_shox}), to the set $K_{1}^{+}(t_{1}:t_{2};\mathbf{r}(t_{0}))$ comprised of the union of single ballistic arcs with $\| \Delta \mathbf{v} \| \le \Delta v_{tot}$. If we stipulate that the spacecraft be capable of maneuvers for which  $\| \Delta \mathbf{v} \| = \Delta v_{tot}$ at the cone vertex, a particularly simple relation holds:
\newline

\noindent Theorem 2: 
\emph{$K^{+}_{\infty}(t_{1}:t_{2};\mathbf{r}(t_{0})) = K_{1}^{+}(t_{1}:t_{2};\mathbf{r}(t_{0}))$}  

 \emph{Proof}:  \emph{$K_{1}^{+} \subset K^{+}_{\infty}$}:  By their respective definitions.
 
 \emph{$K^{+}_{\infty} \subset K_{1}^{+}$}:  By induction.  We first show that, corresponding to any piecewise ballistic trajectory with $n$ shocks $\Delta \mathbf{v}_{i}, i=0 \cdots n-1$, an equivalent single ballistic trajectory exists connecting the cone vertex with any point on the piecewise trajectory subsequent to the final shock.  We then show the velocity change $ \Delta\mathbf{v}_{0}^{n} $ at the vertex of the single ballistic arc obeys
\begin{equation}
\|  \Delta\mathbf{v}_{0}^{n} \| \le \sum_{i=0}^{n-1}\| \Delta\mathbf{v}_{i} \| \le \Delta v_{tot}.
\end{equation}

In the case $n=1$ of a single shock at time $t_{0}$, it is immediate that
\begin{equation}
\| \Delta\mathbf{v}_{0} \|  \equiv \| \Delta\mathbf{v}_{0}^{0} \| \le \Delta v_{tot}. \label{eq:case1}
\end{equation}
Assume that after $n$ shocks the spacecraft location $\mathbf{r}_{n}(t_{n})$ is connected to the vertex by a single ballistic arc with 
 $\sum_{i=0}^{n-1} \| \Delta \mathbf{v}_{i} \| < \Delta v_{tot}$. We show that after a further shock obeying  $\sum_{i=0}^{n} \| \Delta \mathbf{v}_{i} \| \le \Delta v_{tot}$, any point on the
resulting ballistic orbit of the spacecraft  $\mathbf{r}_{n+1}(t)$  with $t > t_{n+1}$ is likewise connected to the cone vertex by a single ballistic arc $\in K_{1}^{+}(t_{1}:t_{2};\mathbf{r}(t_{0}))$.

Take any point $\mathbf{r}_{n+1}(t)$ in $(t_{n+1},t_{2})$. We seek a single ballistic arc connecting 
the cone vertex and $\mathbf{r}_{n+1}(t),\mathbf{v}_{n+1}(t)$ of the form
\begin{equation}
\left ( \begin{array}{c}
\mathbf{r}_{n+1}(t) \\  \mathbf{v}_{n+1}(t)  \label{eq:case_n+1}
\end{array} \right ) =
 \mathbf{\Phi}^{n+1} 
\left ( \begin{array}{c}
\mathbf{r}_{0} \\  \mathbf{v}_{0}^{n+1}
\end{array} \right )
\end{equation}
with 
 \begin{equation}
 \mathbf{v}^{n+1}_{0}  \equiv   \mathbf{v}_{0}+ \Delta \mathbf{v}_{0}^{n+1} 
 \end{equation}
for Lagrangian coefficients $\mathbf{\Phi}^{n+1}$ corresponding to 
some ballistic orbit
$\subset K^{+}_{1}(t_{1}:t_{2};\mathbf{r}(t_{0}))$.
To this relation we append the condition that determines the true anomaly corresponding to the final time of the single ballistic arc,
\begin{equation}
t(\theta)-t_{0}=\sum_{j=1}^{n+1} (t_{j}-t_{j-1})=(t_{n+1}-t_{n})+\sum_{j=1}^{n} (t_{j}-t_{j-1}) 
\label{eq:r0nc}.
\end{equation}
Equations (\ref{eq:case_n+1}) and (\ref{eq:r0nc}) comprise
(in light of the dependence of the Lagrange coefficients upon $\mathbf{v}^{n+1}_{0}$)
a nonlinear vector equation in four unknowns of the form
\begin{equation}
\zeta(\theta,\mathbf{v}_{0}^{n+1},t,\mathbf{r}_{n+1}, \mathbf{v}_{n+1})=0. \label{eq:zetaEQ0}
\end{equation}
The relation (\ref{eq:zetaEQ0}) between $(\theta,\mathbf{v}_{0}^{n+1})$ and $(t,\mathbf{r}_{n+1}, \mathbf{v}_{n+1})$ is, with stipulated $\mathbf{r}_{0}$, diffeomorphic.  As a result, a 
version of the implicit function theorem proved by Kumagai~\cite{Kumagai1980}  
supplies a (real) solution of (\ref{eq:zetaEQ0}) for $(\theta,\mathbf{v}_{0}^{n+1})$.
 Thus, there exists a single ballistic arc corresponding to the $n+1$ shock piecewise ballistic trajectory.

It remains to estimate the magnitude of 
\begin{equation}
\Delta\mathbf{v}_{0}^{n+1} = \mathbf{v}_{0}^{n+1}-\mathbf{v}_{0}.
\end{equation}
The inductive hypothesis is that a single ballistic arc exists corresponding to a trajectory which has experienced $n$ shocks, for which 
\begin{equation}
\sum_{i=0}^{n-1}\| \Delta\mathbf{v}_{i} \| <  \Delta v_{tot}.
\end{equation}
and
\begin{equation}
 \| \Delta\mathbf{v}_{0}^{n} \| \le \sum_{i=0}^{n-1}\| \Delta\mathbf{v}_{i} \|. \label{eq:casen}
\end{equation}
Assume that (\ref{eq:casen}) holds for $n > 1$. Now consider the full trajectory with $n+1$ shocks, subject to the overall limit 
\begin{equation}
\sum_{i=0}^{n}\| \Delta\mathbf{v}_{i} \| \le  \Delta v_{tot}.
\end{equation}
We may replace that portion of the trajectory corresponding to the first $n$ shocks by the equivalent single ballistic arc. The entire trajectory may thus be replaced by one comprised of two ballistic arcs.  Applying  (\ref{eq:case1}) and (\ref{eq:casen}), we have
\begin{eqnarray}
\|  \Delta\mathbf{v}_{0}^{n+1} \| \le \| \Delta\mathbf{v}_{0}^{n} \| + \|  \Delta\mathbf{v}_{n} \| \\
\le \sum_{i=0}^{n-1}\| \Delta\mathbf{v}_{i} \| + \|  \Delta\mathbf{v}_{n} \|
\end{eqnarray}
We conclude 
\begin{equation}
\|  \Delta\mathbf{v}_{0}^{n+1} \| \le \sum_{i=0}^{n}\| \Delta\mathbf{v}_{i} \| \le \Delta v_{tot}.
\end{equation}

Thus, a single ballistic arc  
$\subset K_{1}^{+}(t_{1}:t_{2};\mathbf{r}(t_{0}))$ 
 connects each
 $\mathbf{r}_{n} \in K_{n}^{+}(t_{1}:t_{2};\mathbf{r}(t_{0}))$ with the cone vertex.   
 Call the set of
 all such single ballistic arcs $K_{1}^{+*}$.  
 Therefore, 
\begin{eqnarray}
K^{+}(t_{1}:t_{2};\mathbf{r}(t_{0})) = \bigcup_{n=1}^{\infty} K_{n}^{+}(t_{1}:t_{2};\mathbf{r}(t_{0})) \\
\subset  K_{1}^{+*}  \subset K_{1}^{+},
\end{eqnarray}
proving the theorem. $\,\square$

$K^{+}_{\infty}$ is a subset of $B([t_{0},t_{1}])$, the set of bounded functions on $[t_{0},t_{1}]$, which is a complete metric space in the sup norm. $K^{+}_{\infty}$ is also \emph{closed}:  The mapping from 
\begin{equation}
\left ( \begin{array}{c}
\mathbf{r}_{0} \\  \mathbf{v}_{0}^{n}
\end{array} \right )
\end{equation}
to
\begin{equation}
\left ( \begin{array}{c}
\mathbf{r}_{n}(t) \\  \mathbf{v}_{n}(t)  \label{eq:mapping}
\end{array} \right ) 
\end{equation}
in $K^{+}_{1} \equiv K^{+}_{\infty}$ is homeomorphic (diffeomorphic, in fact) and thus maps closed sets onto closed sets.  The  set  $\{\mathbf{v}_{0}^{n} \}$ that generates $K^{+}_{1}$ is delimited by $ \Delta  \mathbf{v}_{0}^{n} \le \| \Delta v_{tot} \|$; it, and thus $K^{+}_{\infty}$, is a closed set. The latter is therefore a complete metric space.  

Theorem 2 gives us the future cone $K^{+}_{\infty}$ for a countable number of shocks in terms of the trajectories resulting from a single shock experienced at 
$\mathbf{r}(t_{0})=\mathbf{r}_{0}$.  The generic multishock trajectory amounts to a Devil's staircase in $\Delta v$ proceeding from $\Delta v=0$ to $\Delta v=\Delta v_{final} \le \Delta v_{tot}$. The theorem to be proved next extends this result to the case of
a spacecraft maneuvering by continuous nongravitational acceleration as, for instance, in the case of the powered 
flight of a rocket: 
\newline

\noindent Theorem 3: 
\emph{Let $\mathbf{r}(t), t  \in [t_{0},t_{1}]$ the trajectory of a continuously accelerating spacecraft, originating at $\mathbf{r}_{0}$ with initial velocity $\mathbf{v}_{0}$, and subject to $\Delta \mathbf{v} \le \Delta \, v_{tot}$. Then $\mathbf{r}(t)$  $\in K^{+}_{\infty}(t_{1}:t_{2};\mathbf{r}(t_{0}))$}

 \emph{Proof}: Let $\mathbf{r}(t)$ be the trajectory, originating at $\mathbf{r}(t_{0})$ with velocity $\mathbf{v}_{0}$, that results from an imposed acceleration
 \begin{equation}
 \mathbf{a}=\frac{ \mathbf{F}}{m}.
 \end{equation}
 We assume that $ \mathbf{a}$ is a 
 integrable function of $t$.
 The motion of the spacecraft is given by  (\ref{eq:eom0}):
\begin{eqnarray}
 \frac{d}{dt} \left ( \begin{array}{c}
 \mathbf{r} \\ \mathbf{v} 
  \end{array} \right ) =
 \left ( \begin{array}{c}
 \mathbf{v} \\ -\mu \frac{\mathbf{r}}{r^{3}} +\mathbf{a}(t)
  \end{array} \right ). \label{eq:eom}
\end{eqnarray}
 
The effect of continuous acceleration is approximated by a series of shocks 
 $\Delta \mathbf{v}_{i}$ from which we construct a minimizing sequence of shock histories.~\footnote{
 Compare the use of segments of a classical orbit in constructing the quantum-mechanical path integral, \emph{vide.} p. 34 of \emph{Quantum Mechanics and Path Integrals} (1965), R. Feynman and A. Hibbs, McGraw-Hill, New York. } Let the trajectory resulting from the $n$th choice of shock history be given by
\begin{eqnarray}
 \frac{d}{dt} \left ( \begin{array}{c}
 \mathbf{r}_{n} \\ \mathbf{v}_{n} 
  \end{array} \right ) =
 \left ( \begin{array}{c}
 \mathbf{v}_{n} \\ -\mu \frac{\mathbf{r}_{n}}{r_{n}^{3}} +\mathbf{a}_{n}(t)
  \end{array} \right ) 
  \end{eqnarray}
  with $\mathbf{a}_{n}(t)$ a collection of impulsive velocity changes punctuating ballistic motion 
 subject to gravitation alone.  
 We may approximate  $\mathbf{a}_{n}(t)$  for $t \in [t_{0},t_{1}]$ by a sequence of step functions 
  \begin{equation}
 \mathbf{a}_{n}(t)=\sum_{i=1}^{k_{n}} \chi_{[t_{0},t_{1}]}( \tau_{i}) \mathbf{a}^{i}_{n}
 \end{equation} 
 as $k_{n} \rightarrow \infty$, where $\chi_{A}(x)$  is the characteristic function  for an a subset $x \subset A$ 
 and the support of $\mathbf{a}_{n}^{i}$ is
 $supp(\mathbf{a}_{n}^{i}) \equiv \tau_{i}$.
The acceleration $\mathbf{a}_{n}(t)$ integrates to $\Delta \mathbf{v}(t)$ of the form
\begin{equation}
 \int_{t_{0}}^{t} ds  \,\mathbf{a}_{n}=\sum_{i=1}^{k_{n}} \mu( \tau_{i}) \mathbf{a}^{i}_{n} \label{eq:stepfunction}
\end{equation} 
where $\mu(x)$ is the (upper) measure of the set  $x \subset [t_{0},t_{1}]$.
 Writing
\begin{eqnarray}
 \begin{array}{c}
 \mathbf{r}= \mathbf{r}_{n}+\delta \mathbf{r} \\
 \mathbf{v}= \mathbf{v}_{n}+\delta \mathbf{v} \\
 \mathbf{a}= \mathbf{a}_{n}+\delta \mathbf{a} \\
\end{array} \label{eq;deltadefs}
\end{eqnarray} 
we have
\begin{eqnarray}
 \frac{d}{dt} \left ( \begin{array}{c}
\delta  \mathbf{r} \\\delta  \mathbf{v} 
  \end{array} \right ) =
 \left ( \begin{array}{c}
 \delta \mathbf{v} \\ -\mu (\frac{\mathbf{r}}{r^{3}}-
 \frac{\mathbf{r}-\delta  \mathbf{r}}{|\mathbf{r}-\delta  \mathbf{r}|^{\frac{3}{2}}})
  +\delta \mathbf{a}(t)
  \end{array} \right ). 
\end{eqnarray}
 Expanding to second order in $\frac{\delta  r}{r}$,
\begin{eqnarray}
 \frac{d}{dt} \left ( \begin{array}{c}
\delta  \mathbf{r} \\\delta  \mathbf{v} 
  \end{array} \right ) =
 \left ( \begin{array}{c}
 \delta \mathbf{v} \\ -\mu (\delta  \mathbf{r}+3 \frac{\mathbf{r} \cdot \delta  \mathbf{r}}{r^{2}}
 \mathbf{r}) +\delta \mathbf{a}(t)
  \end{array} \right ) + 
  \left ( \begin{array}{c} 0 \\ \mathbf{F}(t)   \end{array} \right ) \label{eq:linearized}
\end{eqnarray}
with
 \begin{equation}
\mathbf{F}(t)= \frac{\mu}{r^{3}} \left (3 \frac{\mathbf{r} \cdot \delta  \mathbf{r}}{r^{2}} \delta  
\mathbf{r} +\frac{3}{2}\frac{ \delta  \mathbf{r} \cdot \delta  \mathbf{r}}{\mathbf{r} \cdot  \mathbf{r}} 
\mathbf{r}+
\frac{15}{2} \frac{(\mathbf{r} \cdot \delta  \mathbf{r})^{2}}{(\mathbf{r} \cdot  \mathbf{r})^{2}}
\mathbf{r} \right )+ O\left ( (\frac{\delta  r}{r})^{3} \right ) \label{eq:Fdef}
\end{equation}

We note that control laws for the shock history $\Delta \mathbf{v}_{i}$ exist  that drive $\delta  \mathbf{r}$ to small values
 \begin{equation}
 \frac{ \delta  r }{r}  \ll 1.
%\frac{| \delta  \mathbf{r} |}{r}  \ll 1.
\end{equation}
It may thus be assumed without loss of generality that if
  \begin{equation}
 \| \delta \mathbf{r} \|_{\infty} \le \alpha  \| \mathbf{r} \|_{\infty}  \label{eq:alphadef}
 \end{equation}
 in the sup norm on $[t_{0},t_{1}]$ then $\exists \, \alpha > 0$ such that the second-order terms in
(\ref{eq:linearized}) dominate the error resulting from linearization.  Using (\ref{eq:alphadef}) and the Schwartz inequality, we have on $[t_{0},t_{1}]$
\begin{eqnarray}
 \begin{array}{c}
\delta  \mathbf{r} \cdot \delta  \mathbf{r} \le 3 \|\delta  \mathbf{r}  \|^{2}_{\infty} \\
 \le 3 \alpha^{2} \|  \mathbf{r} \|^{2}_{\infty} \\
 and \\
\mathbf{r} \cdot \delta  \mathbf{r}  \le 3 \alpha \|  \mathbf{r} \|^{2}_{\infty}.
\end{array} \label{est90}
\end{eqnarray}
We also have
 \begin{eqnarray}
 \begin{array}{c}
\delta   \mathbf{r} \cdot \delta  \mathbf{r}   \le \sqrt{3}  \|\delta   \mathbf{r} \|_{\infty}  
\sum_{j} | \delta r_{j} | \\
\le \sqrt{3}  \alpha \| \mathbf{r} \|_{\infty}  \sum_{j} | \delta r_{j} | \\
 and \\
\mathbf{r} \cdot \delta  \mathbf{r} \le \sqrt{3}  \| \mathbf{r} \|_{\infty}  \sum_{j} | \delta r_{j} |
\end{array}  \label{est91}
\end{eqnarray}
where the sum over $j$ is a sum over the components of $| \delta \mathbf{r} |$. Then, using the estimates (\ref{est90}) and (\ref{est91}), a reliable 
overestimate of the second-order error that is linear in the components of $| \delta  \mathbf{r} |$ is
\begin{equation}
\mathbf{F}(t) < \frac{\mu}{r^{3}} 
\left (9 \frac{\alpha \|  \mathbf{r}  \|^{2}_{\infty}}{r^{2}} | \delta  \mathbf{r} | +
\frac{3 \sqrt{3}}{2}\frac{ \alpha \| \mathbf{r} \|_{\infty}^{2}  \sum_{j} | \delta r_{j} |} {r^{2}} \mathbf{1} +
\frac{45 \sqrt{3}}{2} \frac{ \alpha \|  \mathbf{r}\|^{4}_{\infty} 
\sum_{j} | \delta r_{j} |}{r^{4}} \mathbf{1} \right ) \label{eq:overF}
\end{equation}
with
\begin{equation}
\mathbf{1} \equiv  \left ( \begin{array}{c}
 1 \\
 1\\
 1
 \end{array} \right ). 
\end{equation} 
Upon substitution from (\ref{eq:overF}) and replacement of
\begin{equation}
-3 \frac{\mathbf{r} \cdot \delta  \mathbf{r}}{r^{2}} \mathbf{r}
 \end{equation}
by
\begin{equation}
3 \frac{ \mathbf{r}  \cdot | \delta  \mathbf{r} |}{r^{2}}  \| \mathbf{r} \|_{\infty} \mathbf{1}, 
 \end{equation}
 (\ref{eq:linearized}) takes the form
 \begin{eqnarray}
 \frac{d}{dt} \left ( \begin{array}{c}
 \delta  \mathbf{r} \\ \delta  \mathbf{v} 
  \end{array} \right ) \le \mathbf{A}(t)
 \left ( \begin{array}{c}
| \delta  \mathbf{r} | \\  \delta  \mathbf{v} 
  \end{array} \right ) +
 \left ( \begin{array}{c}
0 \\ \delta  \mathbf{a}(t) 
 \end{array} \right )   
\end{eqnarray}
which, upon formal integration from $t_{0}$ to $t$ becomes
\begin{eqnarray}
 \left ( \begin{array}{c}
\delta  \mathbf{r} \\ \delta  \mathbf{v} 
  \end{array} \right ) \le \int_{t_{0}}^{t} ds \, \mathbf{A}(s)
 \left ( \begin{array}{c}
| \delta  \mathbf{r} |\\ \delta  \mathbf{v} 
  \end{array} \right ) +
 \left ( \begin{array}{c}
0 \\ \delta  \mathbf{w}(t) 
 \end{array} \right )   \label{eq:Volterra}
\end{eqnarray}
where elements of the matrix $\mathbf{A}(t) \ge 0$ on $[t_{0},t_{1}]$ and we are at liberty to choose
$\mathbf{a}_{n}$ so that
\begin{equation}
 \delta  \mathbf{w}(t)=\int_{t_{0}}^{t} ds  \, \delta \mathbf{a}(t). 
 \end{equation}
 is positive on $[t_{0},t_{1}]$ for sufficiently large $n$.

Using $\mathbf{r}_{n}= \mathbf{r}+\delta \mathbf{\rho}$ in place of 
$\mathbf{r}= \mathbf{r}_{n}+\delta \mathbf{r}$ in (\ref{eq;deltadefs}), we find the expected
changes to (\ref{eq:linearized}):  The terms linear in $\delta \mathbf{\rho}$ change sign compared to  (\ref{eq:linearized}), while the form of
$\mathbf{F}(t)$ remains unaltered.  Following the same development leading to (\ref{eq:Volterra}),
with the change that 
\begin{equation}
3 \frac{ \mathbf{r}  \cdot | \delta  \mathbf{\rho} |}{r^{2}}  \| \mathbf{r} \|_{\infty} \mathbf{1}
 \end{equation}
now replaces
\begin{equation}
3 \frac{\mathbf{r} \cdot \delta  \mathbf{\rho}}{r^{2}} \mathbf{r},
 \end{equation}
 and setting $\delta  \mathbf{\rho} = - \delta  \mathbf{r}$, we obtain
\begin{eqnarray}
 \left ( \begin{array}{c}
-\delta  \mathbf{r} \\ \delta  \mathbf{v} 
  \end{array} \right ) \le \int_{t_{0}}^{t} ds \, \mathbf{A}(s)
 \left ( \begin{array}{c}
| \delta  \mathbf{r} | \\  \delta  \mathbf{v} 
  \end{array} \right ) +
 \left ( \begin{array}{c}
0 \\ \delta  \mathbf{w}(t) 
 \end{array} \right )   \label{eq:Volterra2}
\end{eqnarray}
with the identical form for the positive matrix function $\mathbf{A}(s)$.  Both $\delta  \mathbf{r}$ and
$- \delta  \mathbf{r}$ obey the same vector inequality.  It must be the case that $| \delta  \mathbf{r} |$ also obeys the inequality (\ref{eq:Volterra}). 
 The multivariate version of the Gronwall inequality
 ~\cite{ChandraDavis1975,Walter1970,ChuMetcalfe1967} then gives
\begin{eqnarray}
 \left ( \begin{array}{c}
| \delta  \mathbf{r} | \\ \delta  \mathbf{v} 
  \end{array} \right ) \le \int_{t_{0}}^{t} ds \,\mathbf{V}(s,t) 
 \left ( \begin{array}{c}
0 \\ \delta  \mathbf{w}(s) 
  \end{array} \right ) +
 \left ( \begin{array}{c}
0 \\ \delta  \mathbf{w}(t) 
 \end{array} \right )   \label{eq:Gron}
\end{eqnarray}  
where $\mathbf{V}(s,t)$ satisfies
\begin{equation}
\mathbf{V}(s,t) = \mathbf{I}+  \int_{s}^{t}\mathbf{A}(r) \mathbf{V}(r,s) \, dr .
\end{equation} 
In particular, (\ref{eq:Gron}) gives us
  \begin{equation}
  \| \delta  \mathbf{r} \|_{\infty} \le \left \| \int_{t_{0}}^{t_{1}} ds \, 
    \mathbf{V}(s)  \left ( \begin{array}{c}
0 \\ \delta  \mathbf{w} 
  \end{array} \right ) \right \|_{\infty}. \label{eq:ugly}
 \end{equation} 
 %To show (\ref{eq:Lipschitz}), note that
% \begin{equation} 
%  \mathbf{V}
 %\left ( \begin{array}{c}
%0 \\ \delta  \mathbf{w} 
%  \end{array} \right )= \left [  \mathbf{V} \right ]_{j,k}  \delta  w_{j}
% \end{equation} 
% for $i=1,6$ and $j=4,6$.  Therefore,
% \begin{eqnarray} 
% \begin{array}{c}
% \left \|  \int_{t_{0}}^{t_{1}} ds \, \left [  \mathbf{V}(s)
%\left ( \begin{array}{c}
%0 \\ \delta  \mathbf{w} 
%  \end{array} \right )   \right ]  \right \| = \max_{i} \left \| \int_{t_{0}}^{t_{1}} ds 
%  \, \left [  \mathbf{V}(s) \right ]_{i,j}  \delta  w_{j} \right \| \\
%  =\left \| \sum_{j} \int_{t_{0}}^{t_{1}} ds 
%  \, \left [  \mathbf{V}(s) \right ]_{i,j}  \delta  w_{j} \right \| 
%  \end{array}
% \end{eqnarray}  
% for some $i$ in $1,6$ and $j=4,6$.   We have from the Schwartz inequality
 %\begin{eqnarray}  
%\left \| \sum_{j=4} ^{6} \int_{t_{0}}^{t_{1}} ds 
%  \, \left [  \mathbf{V}(s) \right ]_{i,j}  \delta  w_{j} \right \|  \le \sum_{j=4} ^{6} \left \|  \int_{t_{0}}^{t_{1}} ds 
%  \, \left [  \mathbf{V}(s) \right ]_{i,j}  \delta  w_{j} \right \| \\
%  \le \sum_{j=4} ^{6} \sqrt{\left |  \int_{t_{0}}^{t_{1}} ds 
%  \, \left [  \mathbf{V}(s) \right ]_{i,j}^{2} \right | \left  | \int_{t_{0}}^{t_{1}} ds \, \delta  w_{j}^{2}  \right | }.
% \end{eqnarray}  
%By the mean value theorem for integrals,
%\begin{equation}
%\int_{t_{0}}^{t_{1}} ds \, \delta  w_{j}^{2} = (t_{1}-t_{0}) \, \delta  w_{j}^{2}(t')
%\end{equation} 
% with $t' \in  [t_{0},t_{1}]$.  Thus,
An elementary but ugly calculation using the Schwartz inequality and the mean value theorem for integrals converts (\ref{eq:ugly}) into the Lipschitz condition
%\begin{equation}  
% \sum_{j=4}^{6} \left \|  \int_{t_{0}}^{t_{1}} ds 
 % \, \left [  \mathbf{V}(s) \right ]_{i,j}  \delta  w_{j} \right \|_{\infty} \le
%  \sum_{j=4}^{6} \sqrt{\left |  \int_{t_{0}}^{t_{1}} ds 
% \, \left [  \mathbf{V}(s) \right ]_{i,j}^{2} \right |} \sqrt{(t_{1}-t_{0})} | \delta  w_{j} |,
 %\end{equation}  
 \begin{equation}  
  \| \delta  \mathbf{r} \|_{\infty} \le C  \| \delta  \mathbf{w} \|_{\infty}. \label{eq:Lipschitz}
 \end{equation}
 Examine now $ \| \delta  \mathbf{w} \|_{\infty}$.  We have, from (\ref{eq:stepfunction}), 
\begin{equation}  
\delta  \mathbf{w}(t)=\lim_{k_{n} \rightarrow \infty} \sum_{i=1}^{k_{n}} \mu( \tau_{i}) \mathbf{a}^{i}_{n}-
\int_{t_{0}}^{t} ds  \,\mathbf{a}(s).
 \end{equation}
 But, $\mathbf{a}$ being integrable, we may choose a shock history approximated by a sequence of step functions such that~\cite{AP1981a}
 \begin{equation}  
\lim_{n, k_{n} \rightarrow \infty} \sum_{i=1}^{k_{n}} \mu(\tau_{i}) \mathbf{a}^{i}_{n} \rightarrow \int_{t_{0}}^{t} ds  \,\mathbf{a}.
 \end{equation}
Therefore, $m$ exists such that for any $\delta > 0$ such that 
\begin{equation}  
\| \delta  \mathbf{w}_{n} \|_{\infty} \le \delta, \forall \, n \ge m,  
\end{equation} 
there is an $\epsilon  > 0$ such that
\begin{equation}  
\| \delta  \mathbf{r} \|_{\infty} \le \frac{\epsilon}{2}, \forall \, n \ge m.   \label{eq:dwupperbound}
\end{equation}
 
The set  
 $X=\{\delta  \mathbf{r}, \forall \, n > m \}$ 
is a bounded subset of the set $C([t_{0},t_{1}])$ of continuous functions defined on $[t_{0},t_{1}]$. 
By  (\ref{eq:dwupperbound}), the error
$| \delta  \mathbf{r} |$
is uniformly bounded on $[t_{0},t_{1}]$.  Thus, if we choose a neighborhood $N_{t}$ of any $t  \in [t_{0},t_{1}]$,  for  all $\tau \in N_{t}$
 \begin{equation}
|  \delta \mathbf{r}(\tau)-\delta \mathbf{r}(t) | \le | \delta \mathbf{r}(\tau) | + | \delta \mathbf{r}(t) | < \epsilon
\end{equation}
uniformly on  $[t_{0},t_{1}]$ for all $\delta \mathbf{r} \in X$. 
$X$ is therefore an equicontinuous set.
The Arzel\`{a}-Ascoli Lemma then implies that a uniformly convergent subsequence $ \delta \mathbf{r}_{n_{k}} $ exists  on $[t_{0},t_{1}]$ whose sup norm tends to zero.~\cite{AP1981b}  
 Therefore, as $n_{k} \rightarrow \infty$, 
$ \mathbf{r}_{n_{k}}(t) \rightarrow \mathbf{r}(t)$
 uniformly on  $[t_{0},t_{1}]$.  
 Recall that $ K^{+}_{\infty}$ is a complete metric space in the sup norm, and thus contains the limit of all its convergent sequences. We may now claim
\begin{equation}
 \mathbf{r}(t) \in K^{+}_{\infty},
 \end{equation}
proving the theorem. $\,\square$

Clearly, the future cone of a maneuvering spacecraft is the union of piecewise ballistic trajectories and continuously accelerating ones, which gives us
\newline  

\noindent{Theorem 4}: \emph{$K^{+}(t_{1}:t_{2};\mathbf{r}(t_{0}))=K^{+}_{1}(t_{1}:t_{2};\mathbf{r}(t_{0}))$}

\emph{Proof}:  \emph{$K^{+}_{1} \subset K^{+}$}:  By their respective definitions.
 
\emph{$K^{+} \subset K^{+}_{1} $}:  By Theorems 2 and 3.  If $K^{+}_{cont.}$ is the set $K^{+}= \{\text{continuously accelerating} \, \mathbf{r}(t) \in K^{+} \}$, $K^{+}_{\infty} \cup K^{+}_{cont.}=K^{+}_{\infty}  \subset K_{1}^{+}$. $\,\square$

\noindent With Theorem 1, Theorem 4 gives immediately  a result long known to duck-hunters.
  
 \section{Application to FengYun-1C interception} \label{FY1C}
 
 On or about 22:26 UTC 11 January 2007, a direct-ascent antisatellite weapon launched from the territory of the People's Republic of China intercepted, and destroyed, the Chinese FengYun-1C weather satellite.~\cite{Pardini07,Johnson08}  FengYun-1C was launched from the Taiyuan Satellite Launch Center on 10 May 1999 into a Sun-synchronous 860 km, $98.8^{\circ}$ orbit.  It was designed for a two-year mission duration.  A replacement spacecraft, FengYun-1D, was launched in May 2002, but FengYun-1C remained operational into 2005 and evidently responded to ground commands as late as January 2007.~\cite{Johnson08} 
 
While the instrument suite and data products of the FengYun meteorological instruments are documented in numerous sources, not many descriptions of the spacecraft have appeared. FengYun-1C is variously described as having a mass between 880 and 960 kg.  The FengYun-1 series spacecraft were three-axis stabilized by a combination of reaction wheels and cold nitrogen jets.  A 1999 report from the Foreign Broadcast Information Service cites an account of FengYun-1A appearing in the Journal of Chinese Society of Astronautics that claims the attitude control system carried 12 kg of nitrogen at launch.~\cite{fbis89}  In an account of the debris cloud resulting from the destruction of FengYun-1C, Ref. \onlinecite{Pardini07} remarks that at launch, the spacecraft mass was 958 kg, but at its end-of-life, the mass was 880 kg.  One may interpret this as an upper limit for the reaction mass available for attitude control of 78 kg, again presumably nitrogen gas.~\footnote{An account in \emph{China Today: Space Industry} (1992), Astronautic Publishing House, Beijing reports that on 14 February 1991 the FengYun-1B spacecraft experienced an anomaly in which most of the attitude control reaction mass was expended in an uncontrolled manner.  The spacecraft was left spinning at 10 RPM with substantial nutation. The description of the mass and geometry of FengYun-1A in Ref.  \onlinecite {fbis89} allows one to estimate its principal moments of inertia, so that (in principle) one can estimate a lower limit to the reaction mass expended.  The resulting estimate, of order 1 kg, is uninformative for present purposes.}
 
The application of the results from the preceding section to analysis of the FengYun-1C interception of January 2007 is straightforward.  One estimates the future cone of the target subsequent to some fiducial time, granted credible assumptions about its residual $\Delta v$, and compares this to a model of the future cone of the interceptor for a stipulated engagement scenario. A practical method for computational purposes is supplied by Monte Carlo. By virtue of Theorem 4, it suffices to pick a number of elements of $K^{+}_{1}$ for each player and (so to speak) observe the fall of shot. 

The initial position and velocity of FengYun-1C at the start of the engagement is calculated from a two-line element set obtained from Spacetrack~\cite{SpaceTrak}, using the Aerospace Corporation's Satellite Orbital Analysis Program (SOAP)~\cite{SOAPref} to propogate the FengYun-1C state to the target cone vertex. We use a model of the ASAT developed by G. Forden of MIT.~\cite{Forden08,Forden06} Forden's model assumes that the Chinese direct-ascent ASAT is based on the two-stage solid-fueled DF-21 IRBM with a small interceptor of mass $600 \, kg$ as a third stage. The DF-21 stages are modeled as 1.4 m in diameter, with stage 1 being 5.1 m in length, stage 2 1.8 m in length. A center-perforated grain consisting of a double base propellant  is assumed for both stages, with an $I_{sp}$ of 225 seconds for the first stage and an $I_{sp}$ of 230 seconds for the second. Each stage burns for 36 seconds. The model neglects $\Delta v$ of the third stage.

The vertex for the  ASAT cone was arbitrarily chosen as the point in its flyout at which the ASAT altitude exceeds 104 km, using the same launch azimuth and initial pitch assumed in Forden's analysis.  At this point in the flight profile, the missile is late in its second stage burn.  The total $\Delta v$ used to calculate $K^{+}_{ASAT}$ is obtained from (\ref{eq:rocketeqn}) using the partially expended second plus third stage mass for $m_{i}$ and the  second stage tare plus third stage mass for $m_{f}$.
%an estimate of the third stage interceptor tare mass for $m_{f}$. 
%The contribution to $\Delta v_{tot}$ from the interceptor is presumed to be small, and is neglected. 
%The contribution to $\Delta v_{tot}$ from the interceptor  is neglected. 
This procedure presumably overestimates slightly the $\Delta v$ likely available to the actual ASAT very near the cone vertex, but should be accurate enough for illustrative purposes.
 
The change in rocket velocity $\Delta v$ resulting from the exhaust of reaction mass is given by the rocket equation
\begin{equation}
\Delta v=I_{sp} \,  g_{0} \, log\left (\frac{m_{i}}{m_{f}} \right ), \label{eq:rocketeqn}
 \end{equation}
where $I_{sp}$ is the specific impulse in seconds, $g_{0}$ is the gravitational acceleration at 
sea level, and $m_{i}$ and $m_{f}$ are the initial and final masses of the rocket, respectively. The ideal specific impulse for a cold nitrogen gas jet is given as 76 seconds in Ref. \onlinecite{Sutton92}, p. 229.  The range 12 kg-78 kg for the estimated available reaction mass equates to a total $\Delta v$ in the range $11 \, m/s$ to $63 \, m/s$, assuming a tare mass of 880 kg for FengYun-1C. The smaller value of $11 \, m/s$ is used to calculate the FengYun-1C future cone $K^{+}_{FY1C}$ on the ground that this value corresponds to an overestimate of the likely $\Delta v$ available to the spacecraft. Even were the total $\Delta v=63 \, m/s$ initially, by January 2007 the reaction mass would almost certainly have been quite depleted. The vertex of the FengYun-1C cone is chosen at the time of ASAT launch.

\begin{figure}[p]
\includegraphics[width=180mm]{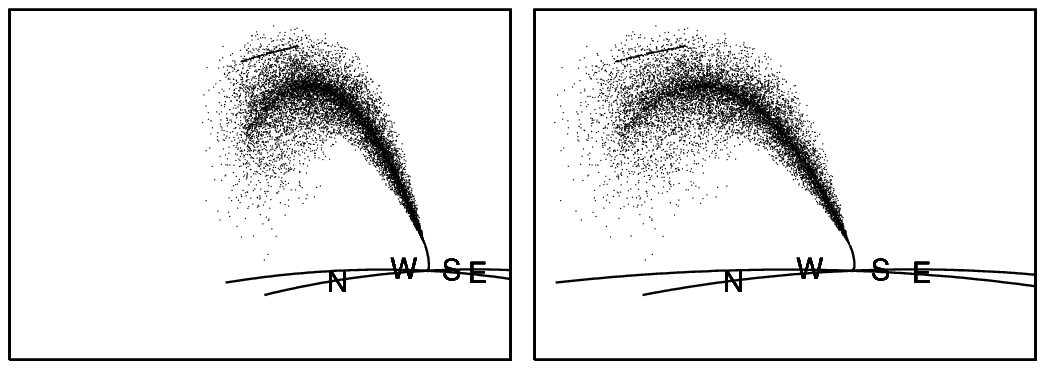}
\caption{Future cones  $K^{+}_{FY1C}(425:475;\mathbf{r}(0))$ and 
$K^{+}_{ASAT}(68:750;\mathbf{r}(68))$ for target and interceptor, respectively, in an Earth-centered inertial frame with times in seconds after interceptor launch. 
The stereogram also shows the latitude and longitude lines of the assumed interceptor launch from $28.13^{\circ} N, 102.02^{\circ} E$ (launch azimuth $345.73^{\circ}$) at 22:18:30 UTC 11 January 2007 and the interceptor trajectory from launch point to $K^{+}_{ASAT}$ vertex. The typical  interceptor cone trajectory points into the page, so that the launch point appears nearer than the interception region.
Distal termini for $K^{+}_{ASAT}$ are truncated at altitudes below 90 km.}  \label{fig1}
\end{figure}

\begin{figure}[p] 
\includegraphics[width=180mm]{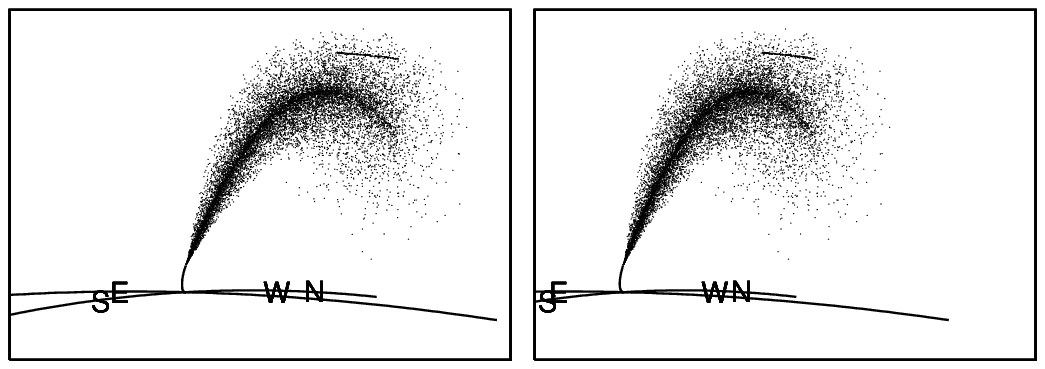}
\caption{Future cones  $K^{+}_{FY-1C}$ and $K^{+}_{ASAT}$ for target and interceptor, respectively, as in Figure  \ref{fig1}. In this view, the  typical interceptor cone trajectory points somewhat out of the page.}  \label{fig2}
\end{figure}

The superposition of the future cone $K^{+}_{ASAT}$ of the ASAT and that part of the cone $K^{+}_{FY1C}$ for FengYun-1C lying within $K^{+}_{ASAT}$ is presented in Figures \ref{fig1} and \ref{fig2} as random-dot stereograms with approximately antipodal viewing geometries. The nominal encounter time 22:26:00 occurs at  450 seconds TALO.  In the interval during which the cones intersect, 
 \begin{equation}
 K^{+}_{FY1C} \subset   K^{+}_{ASAT}
 \end{equation}
 and Theorem 1 guarantees  the ASAT can intercept FengYun-1C. 
 
 That the ASAT could intercept FengYun-1C is hardly a novel conclusion-this much was demonstrated beyond dispute on 11 January 2007.  Admittedly, FengYun-1C appears to have served as a passive target.~\cite{Forden08} Consider, however, a hypothetical scenario in which FengYun-1C maneuvers during the engagement. The disparity in expansion between the cones $K^{+}_{FY-1C}$ and $K^{+}_{ASAT}$  underscores the inability of FengYun-1C to evade the ASAT.  If anything, this exercise understates the vulnerability of FengYun-1C on that date: The cone $K^{+}_{FY1C}$ is computed assuming that the total $\Delta v$ remaining in the attitude control system would actually be available to the spacecraft for maneuver. It is unlikely this assumption holds in practice. At a minimum, the actual  $\Delta v$ is limited by the total thrust the cold nitrogen jets are capable of producing if operated in a configuration that produces net impulse, even on the assumption that operators on the ground issue evasive commands.

\section{Discussion} \label{wottdoesitmean}

 It was noted in Ref. \onlinecite{Morgan2010} that (\ref{eq:Nash2}) in Theorem 1 amounts to establishing a (dominant strategy) Nash equillibrium:  The target may move to any position available in its future, but so long as the interceptor's future contains that of the target, the interceptor can always maneuver to some future position $\mathbf{r}(t_{i})$ of the target.  Cockayne~\cite{ref2}  uses the telling phrase "against all opposition" to describe  the corresponding situation in the game of Two Cars 
 (\emph{vide.}  Isaacs  Ref. \onlinecite{I1965}, p. 202, as well).

Theorem 1 was proved for simply connected future cones, but can be applied to certain multiply connected ones.  If  multiply connected future cones of target and interceptor are manifolds with a timelike foliation and are the finite union of polyhedra of simply connected simplicial complices, one may invoke Theorem 1 for individual polyhedra. Thus, if
\begin{equation}
K^{+}_{int.}(t_{0}:t_{N}; \mathbf{r}^{0}_{p})=\cup_{i=1}^{N-1} K_{int.,i}^{+}(t_{i-1}:t_{i};\mathbf{r}^{0}_{p})
\end{equation}
where each $ K_{int.,i}^{+}$ is a simply connected polyhedron and a subset 
$K^{+}_{targ.}(t_{\alpha}:t_{\omega} ;\mathbf{r}^{0}_{t})$ of $K^{+}_{targ.}$ is a simply connected polyhedron (or other  acyclic absolute neighborhood retract) such that
\begin{equation}
K^{+}_{targ.}(t_{\alpha}:t_{\omega} ;\mathbf{r}^{0}_{t}) \subset K_{int.,i}^{+}(t_{j-1}:t_{ij};\mathbf{r}^{0}_{p})
\end{equation}
for some $j \in 1,N$
then 
the existence of a guaranteed interception opportunity follows in any polyhedron of the interceptor cone that contains $K^{+}_{targ.}(t_{\alpha}:t_{\omega} ;\mathbf{r}^{0}_{t})$. This approach is presumably necessary for a treatment of the Game of Two Cars that relaxes the simplifying conditions regarding initial velocity and cone vertex times made in Section \ref{2cars}.

Theorems 2-4 need not give a good account of the future cone of a spacecraft for which the assumption that it may expend the entirety of its available $\Delta v$ in a single maneuver at the cone vertex is a poor approximation. In that event $K^{+} \subset K^{+}_{1}(\Delta v_{tot})$ remains true, but it may happen that $K^{+} \not \subset K^{+}_{1}(\Delta v_{avail.})$.
A proper treatment of this case with the results of Section \ref{Kepler} requires considering a sequence of future cones, for each of which the $\Delta v_{avail.}$ available in a given time interval of its motion is used in place of $\Delta v_{tot}$. 
%However, Theorems 2-4 should be quite accurate for impulsive maneuvers.
 
\section{Conclusion} \label{byby}

The sample problems worked in Section \ref{futurecones}  show that the future cone construction based on Theorem 1 offers a simple method of determining when a class of differential pursuit/evasion games of kind, in which both players maneuver freely, is guaranteed to have interception of the target as a possible outcome.
 
The method of analysis presented in Section \ref{Kepler} is applicable to direct-ascent and co-orbital ASAT engagements.  The specific example of the FengYun-1C interception by a direct-ascent ASAT demonstrates the applicability of the method to the study of a pursuit/evasion game of considerable practical interest. The formalism developed in this paper and its predecessor may also be taken to apply, if only as a \emph{reducto ad absurdum}, to directed energy weapon attacks with a free line-of-sight to the target satellite.



\begin{thebibliography}{5}

\bibitem{Morgan2010}Morgan, J. A. (2010) Qualitative criterion for interception in a pursuit/evasion game, \emph{Proc. Roy. Soc. A. 466}, pp.1365-1371

\bibitem{Roxin1969}Roxin, E. (1969) Axiomatic Approach in Differential Games, \emph{J. Opt. Theory and App. 3}, 153-163

\bibitem{I1965}Isaacs, R. (1965) \emph{Differential Games}, John Wiley and Sons, New York

\bibitem{F1971}Friedman, A. (1971) \emph{Differential Games}, John Wiley and Sons, New York

\bibitem{ref2}Cockayne, E. (1967) Plane Pursuit with Curvature Constraints,
\emph{SIAM J. Appl. Math. 15}, pp. 1511-1516

\bibitem[Getz and Pachter(1981)]{ref3}Getz, W., and M. Pachter (1981) Capturability in a 
Two-Target "Game of Two Cars", \emph{J. Guidance and Control 4}, pp. 15-21

\bibitem[Miloh(1982)]{ref4}Miloh, T. (1982) A Note on Three-Dimensional Pursuit-Evasion Game with
Bounded Curvature, \emph{IEEE Trans. Automat. Contr. AC-27}, pp. 739-741

\bibitem[Miloh \emph{et. al.}(1993)]{ref5}Miloh T., M. Pachter, A. Segal (1993) The Effect of a
Finite Roll Rate on the Miss-Distance of a Bank-to-Turn Missile, 
\emph{Computers Math. Applic. 26}, pp. 43-53

\bibitem[Shinar and Gutman(1980)]{ref6}Shinar, J., and S. Gutman (1980) Three-Dimensional Optimal
Pursuit and Evasion with Bounded Controls, \emph{IEEE Trans. Automat. Contr. AC-25}, pp. 492-496

\bibitem[Pontryagin(1964)]{P1964}Pontryagin, L. S. (1964) On Some Differential Games, 
\emph{J. SIAM Controls 3}, pp. 49-52 

\bibitem[Pontryagin(1971)]{P1971}Pontryagin, L. S. (1971) Lectures on Differential Games, 
Stanford University, CA (reprinted as DTIC AD724166) 

\bibitem[Pontryagin(1974)]{P1974}Pontryagin, L. S. (1974) On the Evasion Process in Differential 
Games, \emph{Appl. Math and Optimization 1}, pp. 5-19 

\bibitem[Pontryagin(1981)]{P1981}Pontryagin, L. S. (1981) Linear Differential Games of Pursuit, 
\emph{Math. USSR Sbornik 40}, pp. 285-303

\bibitem{EM1946}Eilenberg, S., and D. Montgomery (1946) Fixed Point Theorems for Multi-valued Transformations, \emph{Amer. J. Math 68}, pp.214-222

\bibitem{Tefatsion}Tefatsion, L. (1983) Pure Strategy Nash Equilibrium Points and the Lefschetz Fixed Point Theorem, 
\emph{International Journal of Game Theory 12}, pp. 181-191

\bibitem{K1955}J. L. Kelley (1955) {\it General Topology} (Van Nostrand, Princeton), pp. 270-271.

\bibitem{Cairns46}Cairns, S. (1946) The Triangulation Problem and its Role in Analysis, \emph{Bull. Amer. Math. Soc. 52}, pp. 545-571

\bibitem{Spanier66}Spanier, E., (1966)  \emph{Algebraic Topology}, Springer, New York, pp. 171-172

\bibitem{GranasDj03}Granas, A., and J. Djundji (2003) \emph{Fixed Point Theory}, Springer, New York, pp. 207; 232

\bibitem[Novikov(1967)]{N1967}Novikov, S. P. (1967) Topology of foliations, 
\emph{Trans. Amer. Math. Soc. 14}, pp. 268-304

\bibitem{N1950}J. F. Nash (1950) Equillibrium Points in N-Person Games, {\it Proc.\ Nat.\ Acad.\ Sci.\ USA \bf 36}: 48-49.

\bibitem{N1950a}J. F. Nash (1950) Non-Cooperative Games, dissertation, Princeton University, p. 15

\bibitem[Brooks(2008)]{ref1}Brooks, R. (2008) Game and Information Theory Analysis of Electronic
Countermeasures in Pursuit-Evasion Games, \emph{IEEE Trans. on Syst., Man, Cybern. A 24} pp. 
1281-1294

\bibitem{carter1985}Carter, A (1985) \emph{Anti-Satellite Weapons, Countermeasures, and Arms Control}, U. S. Congress, Office of Technology Assessment OTA-ISC-281, U. S. Government Printing Office, Washington, DC

\bibitem{Battin}Battin, R. (1987) \emph{An Introduction to the Mathematics and Methods of Astrodynamics}, AIAA Educational Series, American Institute of Aeronautics and Astronautics, New York

\bibitem{K1941}Kakutani, S. (1941) A Generalization of Brouwer's Fixed Point 
Theorem, \emph{Duke Math J. 8}, pp. 457-459

\bibitem{Maunder70}Maunder, C. R. F., (1970) \emph{Algebraic Topology}, Cambridge University Press, Cambridge, UK, p. 151

\bibitem{Kumagai1980} Kumagai, S., (1980) An Implicit Function Theorem: Comment, \emph{J. Optimization Th. and App. 31}, pp. 285-288

\bibitem{ChandraDavis1975}Chandra, J., and P. W. Davis (1976) Linear Generalizations of Gronwall's Inequality, \emph{Proceedings of the American Mathematical Society 60}, pp. 157-160

\bibitem{Walter1970}Walter, W. (1970) \emph{Differential and integral inequalities}, Erbebnisse der Mathematik und ihrer Grenzgebiete, Band 55, Springer-Verlag, New York, pp. 143-144

\bibitem{ChuMetcalfe1967}Chu, S. C., and F. T. Metcalfe (1967) On Gronwall's Inequality,
 \emph{Proceedings of the American Mathematical Society 18}, pp. 439-440

\bibitem{AP1981a}Aliprantis, C., and O. Burkinshaw., (1981)  \emph{Principles of Real Analysis}, North Holland, New York,
 p. 153

\bibitem{AP1981b}Aliprantis, C., and O. Burkinshaw., (1981)  \emph{op. cit.}, North Holland, New York,
 pp. 64-65

\bibitem{Pardini07}Pardini, C., and L. Anselmo (2007) Evolution of the Debris Cloud Generated by the FengYun-1C Fragmentation event, \emph{Proceedings of the $20^{th}$ ISSFD}, Annapolis, Maryland

\bibitem{Johnson08}Johnson N., E. Stansbery, J.-C. Liou, M. Horstman, C. Stokely, and D. Whitlock (2008) The characteristics and consequences of the break-up of the FengYun-1C spacecraft, \emph{Acta Astronautica 63}, pp. 128-135

\bibitem{SpaceTrak} http://www.space-track.org

\bibitem{SOAPref}Stodden, D. Y., and G. D. Galasso (1995) Space system visualization and analysis using the Satellite Orbit Analysis Program (SOAP), \emph{IEEE Aerospace Applications Conference Proceedings, 1995}, Aspen, CO, vol 2, pp. 369-387

\bibitem{Forden08}Forden, G. (2008) \emph{A Preliminary Analysis of the Chinese ASAT Test}, unpublished, MIT

\bibitem{Forden06}Forden, G. (2006) $GUI\_Missile\_Flyout$: A General Program for Simulating Ballistic Missiles,  \emph{Science and Global Security, Vol. 15, No. 2},  pp. 133-146

\bibitem{fbis89}Yuhang X. (1989) 'FengYun 1 Meteorological Satellite Detailed, reprinted in 
\emph{China Science and Technology}, JPRS-CST-89-026, Foreign Broadcast Information Service, pp. 1-6 

\bibitem{Sutton92}Sutton, G. (1992) \emph{Rocket Propulsion Elements}, John Wiley and Sons, New York

\bibitem{Marec1979}Marec, J. P. (1979)  \emph{Optimal Space Trajectories}, Elsevier, New York

\end{thebibliography}
\end{document}